\documentclass{article}

\usepackage[english]{babel}

\usepackage{float}
\usepackage{geometry} 
\usepackage{amsmath}
\numberwithin{equation}{section}
\usepackage[ruled]{algorithm2e}
\usepackage{graphicx}
\usepackage[colorlinks=true, allcolors=blue]{hyperref}

\usepackage{xcolor}

\usepackage{authblk}

\usepackage{csquotes}

\usepackage{enumitem}

\usepackage{amsthm}
\usepackage{amssymb}
\usepackage{enumitem}

\newtheorem{theorem}{Theorem}[section]

\newtheorem{lemma}[theorem]{Lemma}

\newtheorem{corollary}{Corollary}[section]

\newtheorem{proposition}[theorem]{Proposition}

\title{A two-point phase recovering from holographic data on a single plane}

\date{}

\author{R.G. Novikov, V.N. Sivkin}

\begin{document}
\maketitle

\begin{abstract}
We consider a plane wave, a radiation solution, and the sum of these solutions
(total solution) for the Helmholtz equation in an exterior region in $\mathbb{R}^d,$ $d\geq 2$.  In this region, we consider a hyperplane $X$ with sufficiently large distance $s$ from the origin in ${\mathbb R}^d.$ We give  two-point local  formulas for approximate recovering the radiation solution restricted to the plane $X$ from the intensity of the total solution at  $X$, that is, from holographic data. The recovering is given in terms of the far-field pattern of the radiation solution with a decaying error term as $s \to +\infty.$ A numerical implementation is also presented.


\end{abstract}


\section{Introduction}

We consider the Helmholtz equation
\begin{align}\label{eq:schrod}
    -\Delta \psi(x) = \kappa^2\psi(x), \quad x \in {\cal U} \subset \mathbb{R}^d, \quad \kappa >0, \quad d\geq 2, 
\end{align}
where ${\cal U}$ is an exterior region in ${\mathbb R}^d,$ that is ${\cal U} = {\mathbb R}^d \setminus (D \cup \partial D),$ where $D$ is an open bounded domain in ${\mathbb R}^d$ with smooth boundary $\partial D,$ and $ {\cal U}$ is connected. An important example is when $D = B_r,$ where 
\begin{align}
    B_r = \{x \in \mathbb{R}^d: |x|<r\}.
\end{align}

Equation \eqref{eq:schrod} arises, in particular, in quantum mechanics, acoustics, and electrodynamics (including optics and X-ray propagation). In addition, equation \eqref{eq:schrod} is a basic equation in holography, going back to \cite{G}.

For equation \eqref{eq:schrod} we consider the plane wave solutions 
\begin{align}\label{eq:psi0}
\psi_0 = e^{ikx}, \quad k \in \mathbb{R}^d, \quad |k| = \kappa, 
\end{align}
 and radiation solutions $\psi_1,$ that is solutions of class $C^2$ which satisfy the Sommerfeld’s radiation condition
\begin{align}\label{eq:sommer}
    |x|^{\frac{d-1}{2}}\left(\frac{\partial}{\partial |x|} - i \kappa\right) \psi_1(x) = 0 \quad \text{ as } |x| \to +\infty,
\end{align}
uniformly in $x/|x|.$ 

Recall that the following asymptotic holds:
\begin{align}\label{eq:1.5}
    \psi_1(x) = \frac{e^{i\kappa|x|}}{|x|^{\frac{d-1}{2}}}f_1\left(\frac{x}{|x|}\right)+ {\cal O}\left(\frac{1}{|x|^{\frac{d+1}{2}}}\right), \quad \text{ as } |x| \to +\infty,
\end{align}
see, for example, \cite{A}, \cite{Karp61}. 

The function $f_1$ arising in \eqref{eq:1.5} is defined on ${\mathbb S}^{d-1}.$ We say  that $f_1$ is the far-field pattern of $\psi_1.$

Let 
\begin{align}\label{eq:X}
    X = X_{s, \omega} = \{ x \in \mathbb{R}^d: \quad (x, \omega) = s\}, \quad s>0, \quad \omega \in \mathbb{S}^{d-1}.
\end{align}

One can see that $X$ is the hyperplane  in $\mathbb{R}^d$  with the normal vector $\omega$ and distance $s$ from the origin $\{0\}$.


In the present work we give  simple local formulas for approximate recovering $\psi_1$ on $X$  from the intensity (hologram) $I =|\psi_0+\psi_1|^2$ on $X.$ 
This recovering is given in term of $f_1$ in \eqref{eq:1.5} with a decaying error term  as $s\to +\infty.$ See Theorem \ref{thm:1}, Propositions \ref{prop:2.3}, \ref{prop:2.5}, \ref{prop:2.2} and Corollaries \ref{cor:2.1}, \ref{cor:2.2} in Section \ref{sec:main}.

By this result, we contribute to holography and phaseless inverse scattering. These studies go back, in particular, to \cite{G}, \cite{Wo1}, \cite{Wo2}. There is a large literature on these problems. In the present work we continue mathematical studies reported, in particular, in  \cite{N1}, \cite{NSh}, \cite{N2023},  \cite{NN}.
In connection with other mathematical works on phaseless inverse scattering and holography, see, e.g., \cite{BD}, \cite{ChS}, \cite{D}, \cite{JL}, \cite{HN},  \cite{HNS}, \cite{K2014}, \cite{K2018}, \cite{KR}, \cite{M}, \cite{MH},   \cite{NS2022}, \cite{Nu}, \cite{R} and references therein.

Recall that the result that $\psi_1$ on $X \subset {\cal U}$ is uniquely determined by the intensity $ I = |\psi_0+\psi_1|^2$ on an arbitrary open domain of $X$ is proved only very recently in \cite{N2023}, \cite{NN} without any approximations. However, determination in \cite{N2023}, \cite{NN} includes an analytic continuation. Therefore, finding more stable reconstructions appropriate for numerical implementations is  a natural objective.
Motivated by this objective, in the present work we continue studies of \cite{NSh}, \cite{N2023}, \cite{NN} by proper adopting  formulas mentioned in Remark 2.5  of \cite{NSh}. The point is that we adopted these  formulas   for approximate local finding $\psi_1$ on $X$ from the intensity $I=|\psi_0+\psi_1|^2$ on $X$ for the case of a single hyperplane $X$, that is, for the basic case in holography.

To our knowledge, approximate  formulas of the present work are  much simpler than related   formulas given in the literature for finding $\psi_1$ on $X$ from the intensity $I = |\psi_0+\psi_1|^2$ on $X.$ In some sense, the formulas of the present work are elementary; in particular, they do not require even any integration. In addition, these formulas do not use any approximation for the Helmholtz equation.

We successfully implemented and illustrated numerically our approximate two-point  reconstruction of $\psi_1$ on $X$ from  $I$ on $X,$ at least, in one of its possible realizations (with respect to the  second point choice) in dimension $d=3;$ see Section \ref{sec:6}. 

The main theoretical results of this work are presented in  details in Section \ref{sec:main} and proved in Sections \ref{sec:proof}, \ref{sec:4}, and  \ref{sec:5}. The aforementioned numerical implementation is presented in Section \ref{sec:6} and also discussed in Section \ref{sec:conc}.

\section{Main results}\label{sec:main}

Let 
\begin{align}\label{eq:8}
    a(x, k) = |x|^{(d-1)/2}(| \psi_0(x)+\psi_1(x)|^2-1), \quad x \in \mathbb{R}^d, \quad k \in \mathbb{R}^d, \quad |k| = \kappa,
\end{align}
where  $\psi_0(x) = e^{ikx}$  as in \eqref{eq:psi0}, $\psi_1$ is a radiation solution for equation \eqref{eq:schrod} as in \eqref{eq:sommer}, \eqref{eq:1.5}.

Let 
\begin{align}\label{eq:sphere_semi}
    {\mathbb S}^+_{\omega} = \{\theta \in \mathbb{S}^{d-1}: \quad (\theta, \omega)>0\}, \quad \omega \in \mathbb{S}^{d-1}.
\end{align}

In particular, we consider $a(x, k)$ at two measurement points $x$ and $y$ on $X = X_{s, \omega},$ where
\begin{align}\label{eq:2.2}
x = x(\theta, s, \omega) = \frac{s\theta}{(\theta, \omega)} \quad \text{and} \quad y = x + \zeta, \quad \zeta \in X_{0, \omega},  \quad \theta \in \mathbb{S}^{+}_{\omega},
\end{align}
 $X_{s, \omega}$ is defined by \eqref{eq:X}.
Note that $\hat{x} = x/|x| = \theta,$ $|x| = s/(\theta, \omega),$ and $x = L_{0, \theta} \cap X_{s, \omega},$ where 
\begin{align}
    L = L_{x_0, \theta} = \{ x \in \mathbb{R}^d: x =  x_0+t \theta: \quad 0< t < +\infty  \}, \quad x_0 \in \mathbb{R}^d, \quad  \theta \in \mathbb{S}^{d-1}.
\end{align}
In addition to $\psi_1,$ we also consider its far-field pattern $f_1$ defined in \eqref{eq:1.5}.

We start with the following result for approximate two-point recovering  $f_1$ on $\mathbb{S}^+_{\omega}$ 
from the intensity $I = |\psi_0+\psi_1|^2$ on the hyperplane $X = X_{s, \omega}$ with sufficiently large $s.$
\begin{theorem}\label{thm:1} Let $\psi_0 = e^{ikx}$ and $\psi_1$ be solutions of equation \eqref{eq:schrod} as in formulas  \eqref{eq:psi0}, \eqref{eq:sommer}, \eqref{eq:1.5}. Let $X = X_{s, \omega}$ be a hyperplane in $\mathbb{R}^d$ as in \eqref{eq:X} with sufficiently large $s$ so that $X \subset {\cal U},$ and $\omega$ is fixed. Then $f_1$ on $\mathbb{S}^+_{\omega}$ is approximately determined by the intensity $|\psi_0 + \psi_1|^2$ on $X$ via 
the following two-point formula: 
\begin{align}\label{eq:2.4}
\begin{aligned}
&f_1(\theta) = f_{1, 1}(\theta) +\frac{1}{D}\left({\cal O}\left(\frac{1}{ |x|^{\sigma}}\right)+{\cal O}\left(\frac{|\zeta|}{ |x|}\right)\right), \text{ as } |x| \to +\infty \text{ and } |x|/|\zeta| \to +\infty, \\
&f_{1, 1}(\theta) := \frac{1}{D}\left(e^{i(ky-|k||y|)}a(x, k)-e^{i(kx-|k||x|)}a(y, k) \right), \\
&D = D(x, \zeta, k)= 2i \sin(k\zeta + |k||x|-|k||y|), \\
&\sigma = 1/2 \text{  for  } d=2, \text{  and  } \sigma = 1 \text{  for  } d\geq 3,
\end{aligned}
\end{align}
where $k$ is fixed,  $a(x, k)$ is defined by \eqref{eq:8}, and $\theta = \hat{x} \in \mathbb{S}^+_{\omega},$ $x, y \in X = X_{s, \omega},$ $\zeta = y - x$ (i.e., $\theta, x, y, \zeta$  are as in \eqref{eq:2.2}), and $\theta$  is fixed. 

\end{theorem}

We consider \eqref{eq:2.4} assuming that $D \neq 0.$ In addition, 
\begin{align}\label{eq:2.5}
    D = 2i \sin \left((k- \kappa\theta) \zeta + {\cal O}\left(\frac{1}{|x|}\right) \right), \text{ as } |x|\to +\infty,
\end{align}
uniformly in 
$\zeta,$ where $|\zeta| \leq \rho,$ where $\rho>0.$ More precisely and more generally, in \eqref{eq:2.4}, 
\begin{align}\label{eq:2.8}
&D  = 
2i \sin\left(  (k -\kappa\theta, \zeta) + \frac{\kappa}{2|x|}\left( (\theta, \zeta)^2 - |\zeta|^2 \right)+{\cal O}\left( \frac{\kappa |\zeta|^3}{|x|^2} \right)\right) , \quad \text{as } \frac{|x|}{|\zeta|}\to+\infty,    
\end{align}
 where $\zeta \in X_{0, \omega}.$

Note that in Theorem \ref{thm:1} and formulas \eqref{eq:2.4}, \eqref{eq:2.8} we have that $|x| = s/(\theta, \omega)\geq s.$ Therefore, formulas \eqref{eq:2.4}, \eqref{eq:2.8} lead to an approximate finding  $f_1$ on $\mathbb{S}^+_{\omega}$ from $|\psi_0+\psi_1|^2$ of $X = X_{s, \omega}$ for sufficiently large $s.$

Theorem \ref{thm:1} and formula \eqref{eq:2.8} are proved in Section \ref{sec:proof}. Formula \eqref{eq:2.5} follows from \eqref{eq:2.8}.

Formula \eqref{eq:1.5}, Theorem \ref{thm:1} and formula \eqref{eq:2.8} imply the following corollary.

 \begin{corollary}\label{cor:2.1} Under the assumptions of Theorem \ref{thm:1}, $\psi_1(x)$ on $X = X_{s, \omega}$ is approximately determined by the intensity  $I = |\psi_0 + \psi_1|^2$ on $X$ via 
the following formula:
\begin{align}\label{eq:2.8_new_0}
    \psi_1(x) = \frac{e^{i\kappa|x|}}{|x|^{\frac{d-1}{2}}}f_{1, 1}\left(\frac{x}{|x|}\right)+\frac{1}{  D} {\cal O}\left(\frac{1}{ |x|^{\frac{d-1}{2}+\sigma} }\right) + {\cal O}\left(\frac{|\zeta|}{|x|^{\frac{d+1}{2}}}\right), \quad x \in X_{s, \omega}, \quad \text{as } |x|\to +\infty \text{ and } \frac{|x|}{|\zeta|} \to +\infty,
\end{align}     
where $f_{1, 1}= f_{1, 1}(\theta),$  $D = D(x, \zeta, k),$ and $\sigma$ are defined in \eqref{eq:2.4}. In addition,
\begin{align}\label{eq:2.9.0}
&D  = 
2i \sin\left(  (k -\kappa\hat{x}, \zeta) + \frac{\kappa}{2|x|}\left( (\hat{x}, \zeta)^2 - |\zeta|^2 \right)+{\cal O}\left( \frac{\kappa |\zeta|^3}{|x|^2} \right)\right) , \quad \text{as } \frac{|x|}{|\zeta|}\to+\infty,
\end{align}
where $\hat{x} = x/|x|,$ $\zeta \in X_{0, \omega}.$ 
 \end{corollary}

One can see that in Theorem  \ref{thm:1} formulas \eqref{eq:2.4} involve just two points $x$ and $y$ of the plane $X_{s, \omega}.$ This explains our terminology "two-point recovering". In addition, this two-point recovering has a property of locality for $\theta \in \mathbb{S}^+_{\omega} \setminus {\cal N}_{\omega, k, \varepsilon},$ where the exceptional domain $ {\cal N}_{\omega, k, \varepsilon}$ is defined by formula \eqref{eq:N_ke} below. This locality means that $ |y-x|$ can be small; see estimate in \eqref{eq:2.9}, where $|\zeta| = |y-x|.$

Let
\begin{align}
&\theta = \theta_{\perp}+ \theta_{||}, \quad k = k_{\perp} + k_{||}, \\   
&{\cal N}_{\omega, k, \varepsilon} = \{\theta\in \mathbb{S}^+_{\omega}: |k_{||}- \kappa \theta_{||}|<\varepsilon\}, \quad k \in \mathbb{R}^d, \quad  |k| = \kappa, \quad \varepsilon \in (0, 2\kappa),  \quad \omega \in \mathbb{S}^{d-1}, 
\label{eq:N_ke} \\
&\mathbb{S}^+_{\omega, \delta} = \{\theta\in \mathbb{S}^+_{\omega}: |\theta_{||}|<\delta\}, \quad \delta \in (0, 1), \label{eq:S+omegadelta} \\
&X_{s, \omega, \delta} = \{x \in X_{s, \omega}: x/|x| \in \mathbb{S}^+_{\omega, \delta}\}, \label{eq:Xsod}
\end{align}

where $k_{\perp} \perp X_{0, \omega},$ $\theta_{\perp}\perp X_{0, \omega},$ and $k_{||}, \theta_{||} \in X_{0, \omega},$ $\omega \in \mathbb{S}^{d-1}.$

Properties of formulas \eqref{eq:2.4} and   \eqref{eq:2.8_new} for approximate finding $f_1$ on $\mathbb{S}^+_{\omega}$ and $\psi_1$ on $X_{s, \omega}$  strongly depend on properties of $D$ in \eqref{eq:2.4}. In particular, we cannot use these formulas for $x$ and $\zeta = y - x,$ if $D = D(x, \zeta, k) = 0$ for fixed $k.$

In connection with properties of $D = D(r\theta, \zeta, k)$ in \eqref{eq:2.4}, we use, in particular, the following Lemmas \ref{lem:2.2} and \ref{lem:2.4}.

\begin{lemma}\label{lem:2.2} Let $\omega \in \mathbb{S}^{d-1},$ $k \in \mathbb{R}^d, $ $|k| = \kappa>0,$ and $\theta \in \mathbb{S}^+_{\omega}.$ Then the following statements hold:

(i) Suppose that $(k, \omega)>0.$  Then $(k-\kappa \theta, \zeta) = 0$ for any $\zeta \in X_{0, \omega},$ if and only if $\kappa \theta =k,$ or equivalently, $\kappa \theta_{||} = k_{||}.$

    (ii) Suppose that $(k, \omega)\leq 0.$ Then $(k-\kappa \theta, \zeta) = 0$ for any $\zeta \in X_{0, \omega},$ if and only if $ \kappa \theta_{||} = k_{||},$ $\kappa \theta_{\perp} = -k_{\perp},$ or equivalently,  just $ \kappa \theta_{||} = k_{||}.$

(iii) In addition, if
\begin{align}\label{eq:2.12}
     |\kappa \theta_{||}-k_{||}|\geq \varepsilon>0, \quad \zeta =   -\alpha\frac{\kappa \theta_{||}-k_{||}}{|\kappa \theta_{||}-k_{||}|^2}, \quad \alpha \in \mathbb{R}\setminus \{0\},
\end{align}
then
\begin{align}
&(k-\kappa \theta, \zeta) = \alpha, \quad  |\zeta|\leq \frac{|\alpha|}{\varepsilon}, \label{eq:2.9} \\
&D(r\theta, \zeta, k) = 2i \sin\left(\alpha +{\cal O}\left(\frac{1}{r}\right)\right), \quad \text{ as } r \to +\infty, \text{ for fixed } \varepsilon \text{ and } \alpha,   \label{eq:2.14.0}  
\end{align}
where $D$ is defined in \eqref{eq:2.4}.

    
\end{lemma}

Note that in \eqref{eq:2.12} we have that $\theta \in \mathbb{S}^+_{\omega} \setminus {\cal N}_{\omega, k, \varepsilon}$ in view of \eqref{eq:N_ke},  and $\zeta \in X_{0, \omega}$  since $\theta_{||},$  $k_{||} \in X_{0, \omega}.$  The purpose of the next lemma consists in choosing 
 $\zeta$ in formulas \eqref{eq:2.4}, \eqref{eq:2.8} for $\theta \in {\cal N}_{\omega, k, \varepsilon}.$ 

\begin{lemma}\label{lem:2.4}
 Let $\omega \in \mathbb{S}^{d-1},$ $k \in \mathbb{R}^d, $ $|k| = \kappa>0,$ $r>0,$ $\theta \in \mathbb{S}^+_{\omega},$ $\alpha <0,$ $\sin \alpha \neq 0.$ Let
 \begin{align}
 &\zeta = -\beta \frac{\kappa \theta_{||}- k_{||}}{|\kappa \theta_{||}- k_{||}|}, \quad  \hat{\zeta} = - sgn \beta \frac{\kappa \theta_{||}- k_{||}}{|\kappa \theta_{||}- k_{||}|},  \text{ if } \kappa \theta_{||}\neq k_{||}, \label{eq:2.14}\\
 &\zeta = \beta \hat{\zeta}, \quad  \hat{\zeta} \in X_{0, \omega}, \quad |\hat{\zeta}| = 1, \quad \text{ if } k\theta_{||} = k_{||}, \label{eq:2.15} \\
 &\beta = \beta(\alpha, \kappa, r, \theta_{||}, k_{||}) =  \frac{ 2 \alpha}{ |k_{||}-\kappa \theta_{||}| +\sqrt{|k_{||}-\kappa \theta_{||}|^2+\frac{2 \kappa}{r} ((\theta_{||}, \hat{\zeta})^2 - 1) \alpha}}. \label{eq:2.18}
 \end{align}
 
    Then the following estimates hold:
 \begin{align}
     &|\zeta| \leq  |\beta| \leq \sqrt{\frac{2\alpha r}{\kappa((\theta_{||}, \hat{\zeta})^2 - 1)}}, \quad \alpha<0,  \label{eq:2.19}\\  
     &D = D(r\theta, \zeta, k) = 2i \sin \left( \alpha + {\cal O}\left( \frac{ 1}{\kappa^{1/2}r^{1/2}}\right) \right), \text{ as } r  \to +\infty,  \label{eq:2.20}
 \end{align}
uniformly in $\theta \in \mathbb{S}^+_{\omega, \delta}$  for fixed $\delta\in (0, 1),$ where $D(x, \zeta, k)$ is defined in \eqref{eq:2.4},  $|x| = r,$ $x/|x| = \theta,$  $\mathbb{S}^+_{\omega, \delta}$ is defined by \eqref{eq:S+omegadelta}. 
\end{lemma}

Note that the estimates of Lemma \ref{lem:2.4} are of interest for $\theta \in {\cal N}_{\omega, k, \varepsilon}.$

Lemmas  \ref{lem:2.2} and \ref{lem:2.4} are proved in Section \ref{sec:4}.

Theorem \ref{thm:1} and Lemmas \ref{lem:2.2} and \ref{lem:2.4} lead to the following propositions for finding $f_1$ on $\mathbb{S}^+_{\omega}$ defined by \eqref{eq:sphere_semi}.

\begin{proposition}\label{prop:2.3}
    Let $\omega \in \mathbb{S}^{d-1},$ $k \in \mathbb{R}^d, $ $|k| = \kappa>0,$ $\alpha \in \mathbb{R},$ $|\sin \alpha| \geq \delta > 0,$ $\varepsilon>0,$ and $\theta \in \mathbb{S}^+_{\omega}.$ Then the following formula holds:
\begin{align}\label{eq:2.13}
&f_1(\theta) = f_{1, 1}(\theta) + \frac{1}{2i \sin \alpha }{\cal O}\left( \frac{1}{|x|^{\sigma}}\right), \quad \text{as } |x| \to +\infty, \quad \theta \in \mathbb{S}^+_{\omega} \setminus {\cal N}_{\omega, k, \varepsilon},
\end{align}
 where $f_{1, 1}(\theta)$ is defined according to \eqref{eq:2.4} and \eqref{eq:2.2} with $\zeta = y - x$ as in   \eqref{eq:2.12}, and ${\cal N}_{\omega, k, \varepsilon}$ is defined by \eqref{eq:N_ke}.  
\end{proposition}

\begin{proposition}\label{prop:2.5}
   Let $\omega \in \mathbb{S}^{d-1},$ $k \in \mathbb{R}^d, $ $|k| = \kappa>0,$  $\theta \in \mathbb{S}^+_{\omega},$ $\alpha <0,$ $\sin \alpha \neq 0.$ Then the following formula holds:
\begin{align}\label{eq:2.14.1}
&f_1(\theta) = f_{1, 1}(\theta) + \frac{1}{2i \sin \alpha }{\cal O}\left( \frac{1}{|x|^{1/2}}\right), \quad \text{as } |x| \to +\infty, \quad \theta \in \mathbb{S}^+_{\omega}, 
\end{align}
 where $f_{1, 1}(\theta)$ is defined according to \eqref{eq:2.4} and \eqref{eq:2.2} with $\zeta = y -x$ as in   \eqref{eq:2.14}, \eqref{eq:2.15}, \eqref{eq:2.18}. In addition, estimate \eqref{eq:2.14.1} holds uniformly in $\theta \in \mathbb{S}^+_{\omega, \delta}$ for fixed $\delta \in (0, 1),$ where $ \mathbb{S}^+_{\omega, \delta}$ is defined by \eqref{eq:S+omegadelta}.
      
\end{proposition}

One can see that $|\zeta| < |\alpha|/\varepsilon$ and is bounded in Proposition \ref{prop:2.3},  whereas $|\zeta|  = {\cal O}(|x|^{1/2})$ in Proposition \ref{prop:2.5}. 

Propositions \ref{prop:2.3} and \ref{prop:2.5} are proved in Section \ref{sec:4}.

Formula \eqref{eq:1.5}, Theorem \ref{thm:1}, and Propositions \ref{prop:2.3} and \ref{prop:2.5} imply the following corollary.

\begin{corollary}\label{cor:2.2}

Under the assumptions of Theorem \ref{thm:1}, $\psi_1(x)$ on $X = X_{s, \omega}$ is approximately determined by the intensity  $I = |\psi_0 + \psi_1|^2$ on $X$ via 
the following formulas:
\begin{align}\label{eq:2.8_new}
&\psi_1(x) = \frac{e^{i\kappa|x|}}{|x|^{\frac{d-1}{2}}}f_{1, 1}\left(\hat{x}\right)+\frac{1}{  2i \sin \alpha} {\cal O}\left(\frac{1}{ |x|^{\frac{d-1}{2}+\sigma} }\right),  \quad x \in X,  \quad \hat{x} = \frac{x}{|x|} \in \mathbb{S}^+_{\omega}\setminus {\cal N}_{\omega, k, \varepsilon}, \quad \text{as } |x| \to +\infty,
\end{align}     
where $ f_{1, 1}(\hat{x}),$  is defined according to \eqref{eq:2.4} and \eqref{eq:2.2} with $\zeta = y - x$ as in   \eqref{eq:2.12}, for $|\sin \alpha|>\delta,$ 
${\cal N}_{\omega, k, \varepsilon}$ is defined by \eqref{eq:N_ke}, $\sigma$ is as in \eqref{eq:2.4},  and $\varepsilon>0$ and $\delta>0$ are fixed;
\begin{align}\label{eq:2.25}
    &\psi_1(x) = \frac{e^{i\kappa|x|}}{|x|^{\frac{d-1}{2}}}f_{1, 1}\left(\hat{x}\right)+\frac{1}{  2i \sin \alpha} {\cal O}\left(\frac{1}{ |x|^{\frac{d}{2}} }\right), \quad x \in X, \quad \hat{x} = \frac{x}{|x|} \in  \mathbb{S}^+_{\omega}, \quad  \text{as } |x| \to +\infty,  
\end{align}
 where $f_{1, 1}(\hat{x})$ is defined according to \eqref{eq:2.4} and \eqref{eq:2.2} with $\zeta = y -x$ as in   \eqref{eq:2.14}, \eqref{eq:2.15}, \eqref{eq:2.18}, for $\alpha<0,$ $\sin \alpha \neq 0.$ In addition, estimate \eqref{eq:2.25} holds uniformly in $\hat{x} \in \mathbb{S}^+_{\omega, \delta}$ for fixed $\delta \in (0, 1),$ where $ \mathbb{S}^+_{\omega, \delta}$ is defined by \eqref{eq:S+omegadelta}.

\end{corollary}

In order to improve the remainder ${\cal O}(|x|^{-1/2})$ in formula \eqref{eq:2.4} for $d = 2,$ we also give
the following result.

\begin{proposition}\label{prop:2.2}
    Suppose that the assumptions of Theorem \ref{thm:1} for $d = 2$ hold.
Then $\psi_1$ on $X$ is approximately determined by  the intensity $|\psi_0 + \psi_1|^2$ on $X$ via formula \eqref{eq:1.5} and the following formula for $f_1:$
\begin{align}\label{eq:2.6}
  f_1(\theta) =  f_{1, 1} - \frac{1}{D}\left(e^{i(ky-|k||y|)}-e^{i(kx - |k||x|)} \right) \frac{|f_{1, 1}|^2}{|x|^{1/2}}   +\frac{1}{D}{\cal O}\left(\frac{|\zeta|}{|x|}\right)+ \frac{1}{D^3}\left({\cal O}\left(\frac{1}{|x|}\right)+{\cal O}\left(\frac{|\zeta|}{|x|^{3/2}}\right)\right), 
\end{align}
as $|x|\to+\infty,$ $|x|/|\zeta|\to +\infty.$

Here, $f_{1, 1} = f_{1, 1}(\theta)$ is defined in \eqref{eq:2.4}, $k$ is fixed,  $x, y,  \theta, \zeta$ are as in \eqref{eq:2.2}, $d=2$, and $\theta$ is fixed.

    
\end{proposition}

Proposition \ref{prop:2.2} is proved in Section \ref{sec:5}.

Formula \eqref{eq:2.6} is a modification of formula \eqref{eq:2.4} for $d=2.$ The reminder ${\cal O}(|x|^{-1})+{\cal O}(|x|^{-3/2}|\zeta|)$ in \eqref{eq:2.6} is better than the reminder ${\cal O}(|x|^{-1/2})$ in \eqref{eq:2.4} for $d=2$. Using formula \eqref{eq:2.6}, one can improve formulas 
\eqref{eq:2.8_new_0}, \eqref{eq:2.13},  \eqref{eq:2.8_new} for $d=2$ in a similar way. 

\section{Proof of Theorem \ref{thm:1} and formula \eqref{eq:2.8}}\label{sec:proof}

Due to \eqref{eq:psi0}, \eqref{eq:1.5}, \eqref{eq:8}, we have that
\begin{align}\label{eq:3.1}
&|\psi_0(x)+\psi_1(x)|^2 = 1 + e^{-ikx}\frac{e^{i\kappa |x|} }{|x|^{\frac{d-1}{2}}} f_1 + e^{ikx}\frac{e^{-i\kappa |x|} }{|x|^{\frac{d-1}{2}}} \overline{f_1} + \frac{1}{|x|^{d-1}} f_1 \overline{f_1} + {\cal O}\left(\frac{1}{|x|^{\frac{d+1}{2}}}\right), \quad \text{as } |x|\to +\infty, \\
&a(x, k) = e^{-ikx+i\kappa |x|} f_1 + e^{ikx-i\kappa |x|} \overline{f_1} + \frac{1}{|x|^{\frac{d-1}{2}}} f_1 \overline{f_1} + {\cal O}\left(\frac{1}{|x|}\right), \quad \text{as } |x|\to +\infty, \label{eq:3.2.0}
\end{align}
where $f_1 = f_1(\theta),$ $\theta = \hat{x} = x/|x|, $ $x\in \mathbb{R}^d,$ uniformly in $\theta \in {\mathbb S}^{d-1}.$

Recall that $f_1$ is smooth on ${\mathbb S}^{d-1};$ see, for example, \cite{A}, \cite{Karp61}.

Proceeding from \eqref{eq:3.2.0} and the aforementioned smoothness of $f_1,$ we obtain the following system of equations for approximate finding $f_1(\theta)$ from $a(x, k)$ and $a(y, k):$ 
 \begin{align}
&\begin{aligned}\label{eq:3.5}
&e^{-ikx}e^{i\kappa |x|} f_1(\theta) + e^{ikx}e^{-i\kappa |x|} \overline{f_1}(\theta) = a(x, k) + {\cal O}\left(\frac{1}{|x|^{\sigma}}\right),   
\end{aligned}\\
&\begin{aligned}\label{eq:3.6.0}
&e^{-iky}e^{i\kappa |y|} f_1(\theta) + e^{iky}e^{-i\kappa |y|}\overline{f_1}(\theta) = a(y, k) + {\cal O}\left(\frac{1}{|x|^{\sigma}}\right) + {\cal O}\left(\frac{|\zeta|}{|x|} \right),
\end{aligned}
\end{align}
 as $|x| \to +\infty$ and  $|x|/|\zeta|\to +\infty,$ uniformly in $\theta\in {\mathbb S}^+_{\omega}.$ Here, $x,$ $y$ are defined as in \eqref{eq:2.2}, and, in particular, $|x| = s/(\theta, \omega);$ $\sigma$ is as in \eqref{eq:2.4}. 

Theorem \ref{thm:1} follows from the linear system \eqref{eq:3.5}, \eqref{eq:3.6.0} for $f_1,$ $\overline{f_1}.$  In particular, $D$ in \eqref{eq:2.4} is the determinant of the  corresponding $2\times 2$ matrix, arising in this system; and it is  assumed that $D\neq 0.$ 

Note that equation \eqref{eq:3.5} follows directly from  \eqref{eq:3.2.0}. In turn, equation \eqref{eq:3.6.0} follows from the formulas: 
\begin{align}
&e^{-iky}e^{i\kappa |y|} f_1(\hat{y}) + e^{iky}e^{-i\kappa |y|} \overline{f_1}(\hat{y}) = a(y, k) + {\cal O}\left(\frac{1}{|y|^{\sigma}}\right),  \\
&f_1(\hat{y}) = f_1(\theta) + {\cal O}\left(\frac{|\zeta|}{|x|}\right), \label{eq:3.6} \\
&\begin{aligned}\label{eq:3.7}
&{\cal O}\left(\frac{1}{|y|^{\sigma}}\right) = {\cal O}\left(\frac{1}{|x|^{\sigma}}\right),
\end{aligned}
\end{align}
as $|x| \to +\infty$ and $|x|/|\zeta| \to +\infty,$ uniformly in $\theta\in {\mathbb S}^+_{\omega}.$ 

Here, formula \eqref{eq:3.6} follows from the aforementioned smoothness of $f$ and the formula 
\begin{align}\label{eq:3.8.1}
\hat{y} = \theta + {\cal O}\left(\frac{|\zeta|}{|x|}\right), \quad \text{as } |x|/|\zeta| \to +\infty,  
\end{align}
formula \eqref{eq:3.7} follows from the formula
\begin{align}\label{eq:3.9.1}
    |y| = |x|\left(1+{\cal O}\left(\frac{|\zeta|}{|x|}\right)\right), \quad  \text{as } |x|/|\zeta| \to +\infty.
\end{align}

This completes the proof of Theorem \ref{thm:1} taking into account that more precise versions of formulas \eqref{eq:3.8.1}, \eqref{eq:3.9.1} are given below by formulas \eqref{eq:3.10}, \eqref{eq:3.9}.

One can see that 
\begin{align}\label{eq:3.10}
&\begin{aligned}
&\hat{y}  = \frac{x+\zeta}{|x+\zeta|} = \frac{|x|\theta+\zeta}{||x|\theta+\zeta|} =  \frac{\theta+\zeta/|x|}{|\theta+\zeta/|x||} = \frac{\theta+\zeta/|x|}{(1+2\frac{(\zeta, \theta)}{|x|} + \frac{|\zeta|^2}{|x|^2})^{1/2}} =  \\
&=(\theta+\zeta /|x|)(1- \frac{(\zeta, \theta)}{|x|}+{\cal O}(\frac{|\zeta|^2}{|x|^2})) = \theta + \frac{\zeta-\theta(\theta, \zeta)}{|x|}+{\cal O}\left(\frac{|\zeta|^2}{|x|^2}\right),    
\end{aligned}\\  
&\begin{aligned}\label{eq:3.9}
&|y| = |x+\zeta| = |x||\theta+\zeta/|x|| = |x|(1+2\frac{(\theta, \zeta)}{|x|} +|\zeta|^2/|x|^2 )^{1/2} = \\
&=|x|\left(1+\frac{(\theta, \zeta)}{|x|} + \frac{|\zeta|^2-(\theta, \zeta)^2}{2|x|^2}  + {\cal O}\left(\frac{|\zeta|^3}{|x|^3}\right)\right),
\end{aligned}
\end{align}
where $|x|/|\zeta|\to +\infty.$ Here, we use that 
$\sqrt{1+\varepsilon} = 1 + \varepsilon/2 - \varepsilon^2/8 + {\cal O}(\varepsilon^3), $ as  $ \varepsilon\to 0. $

Next, formula \eqref{eq:2.8} follows from the formula for $D$ in \eqref{eq:2.4} and the following asymptotic formula:
\begin{align}\label{eq:3.8}
\begin{aligned}
&k\zeta + \kappa(|x|-|y|) = 
k\zeta - \kappa \left((\theta, \zeta) -\frac{|\zeta|^2}{2|x|} + \frac{(\theta, \zeta)^2}{2|x|} +{\cal O}\left( \frac{ |\zeta|^3}{|x|^2} \right)\right)  = \\
&=(k-\kappa \theta, \zeta) +\frac{\kappa}{2|x|}((\theta, \zeta)^2-|\zeta|^2) +{\cal O}\left( \frac{ \kappa |\zeta|^3}{|x|^2} \right), \quad \text{as } |x|/|\zeta|\to+\infty.    
\end{aligned}    
\end{align}
In turn, in \eqref{eq:3.8} we used formula \eqref{eq:3.9} for $|y|$.

\section{Proof of Lemmas \ref{lem:2.2}, \ref{lem:2.4} and Proposition \ref{prop:2.3} }\label{sec:4}

\subsection{Proof of Lemma \ref{lem:2.2}}

In view of our assumptions, we have that 
\begin{align}
    &k - \kappa \theta = k_{||} - \kappa \theta_{||} + k_{\perp} - \kappa \theta_{\perp}, \\
    &k_{||} - \kappa \theta_{||} \in X_{0, \omega}, \label{eq:4.2.1} \\
    &\hat{k}_{\perp} = k_{\perp}/|k_{\perp}| = \pm \omega, \quad \hat{\theta}_{\perp} = \theta_{\perp}/|\theta_{\perp}| =  \omega,  \label{eq:4.3.1}\\
    &k^2_{||}+k^2_{\perp} = \kappa^2, \quad \theta^2_{||} + \theta^2_{\perp} = 1. \label{eq:4.4.2}
\end{align}
Therefore, 
\begin{align}\label{eq:4.4.1}
    &(k-\kappa \theta, \zeta) = (k_{||} - \kappa \theta_{||}, \zeta), \quad \text{for } \zeta \in X_{0, \omega}. 
\end{align}
Due to \eqref{eq:4.2.1}, \eqref{eq:4.4.1}, we have that 
\begin{align}\label{eq:4.6}
    &(k-\kappa \theta, \zeta) = 0, \quad \text{for all } \zeta \in X_{0, \omega}, \text{ if and only if } k_{||} - \kappa \theta_{||} = 0.
\end{align}


In addition, in view of \eqref{eq:4.4.2}, if $k_{||} = \kappa \theta_{||},$ then 
\begin{align}\label{eq:4.7}
\kappa^2 \theta^2_{\perp} = \kappa^2-\kappa^2\theta^2_{||} =\kappa^2 - k^2_{||} = k^2_{\perp}.
\end{align}

Items (i) and (ii) of Lemma \ref{lem:2.2} follow from \eqref{eq:4.3.1}, \eqref{eq:4.6}, \eqref{eq:4.7}. 

Next, for $\zeta$ defined in \eqref{eq:2.12},  the  equality $(k-\kappa \theta, \zeta) = \alpha$ in \eqref{eq:2.9} follows from  \eqref{eq:4.2.1}, \eqref{eq:4.4.1}, and the inequality $|\zeta|\leq \alpha/\varepsilon$ in \eqref{eq:2.9} follows from the inequality $|\zeta| \leq \alpha/|\kappa \theta_{||}-k_{||}|,$ and the assumption that $|\kappa \theta_{||} -k_{||}|\geq \varepsilon.$ 

Finally, due to \eqref{eq:2.9.0} and  \eqref{eq:2.9}, we have that:
\begin{align}
    &D = D(x, \zeta, k) = 2i \sin\left( \alpha + \frac{\alpha^2}{\varepsilon^2|x| } {\cal O}\left(1 \right)\right), \text{ as } \varepsilon |x|/\alpha \to +\infty, \text{ if } \hat{x} = \theta \in \mathbb{S}^+_{\omega} \setminus {\cal N}_{\omega, k, \varepsilon}, \\
    &D = 2i \sin \alpha \left(1 + {\cal O}\left( \frac{\alpha}{\varepsilon^2 |x|}\right)\right), \text{ as }  \varepsilon^2 |x|/\alpha \to +\infty, \quad  \text{ if }  \hat{x} = \theta \in \mathbb{S}^+_{\omega} \setminus {\cal N}_{\omega, k, \varepsilon}. \label{eq:4.9.0}
\end{align}

Formula \eqref{eq:4.9.0} implies \eqref{eq:2.14.0}.

Thus, item (iii) of Lemma \ref{lem:2.2} is also proved.

\subsection{Proof of Lemma \ref{lem:2.4}}

Using formula \eqref{eq:3.8},  for $ \zeta$  in \eqref{eq:2.14}, we obtain that 
\begin{align}\label{eq:4.10.0}
    k \zeta +\kappa(|x|-|y|) = \beta |k_{||}-\kappa \theta_{||}| + \frac{\kappa}{2|x|}\beta^2((\theta, \hat{\zeta})^2 - 1) +{\cal O}\left(\frac{\kappa \beta^3}{|x|^2}\right),
\end{align}
where $x = r \theta.$

We are looking for $\beta$ such that
\begin{align}\label{eq:4.11}
    \beta |k_{||}-\kappa \theta_{||}| + \frac{\kappa}{2r}\beta^2((\theta, \hat{\zeta})^2 - 1) = \alpha.
\end{align}
Using that $(\theta, \hat{\zeta}) = (\theta_{||}, \hat{\zeta}),$ $|\hat{\zeta}| = 1,$ and $|\theta_{||}|<1$ since $\theta \in \mathbb{S}^+_{\omega}, $ we have that 
\begin{align}\label{eq:4.12}
(\theta, \hat{\zeta})^2-1=(\theta_{||}, \hat{\zeta})^2-1<0.    
\end{align}
Due to \eqref{eq:4.12}, equation \eqref{eq:4.11} for $\beta$ can be rewritten in the standard form:
\begin{align}
    \beta^2 + \frac{2r |k_{||}-\kappa \theta_{||}|}{\kappa((\theta_{||}, \hat{\zeta})^2 - 1)}\beta - \frac{2r}{\kappa((\theta_{||}, \hat{\zeta})^2 - 1)} \alpha = 0
\end{align}

 The solution of this quadratic equation is given by 
\begin{align}\label{eq:4.14}
\begin{aligned}
&\beta  = - \frac{r|k_{||}-\kappa \theta_{||}|}{\kappa((\theta_{||}, \hat{\zeta})^2 - 1)} +\sqrt{\frac{r^2|k_{||}-\kappa \theta_{||}|^2}{\kappa^2((\theta_{||}, \hat{\zeta})^2 - 1)^2}+\frac{2r\alpha}{\kappa((\theta_{||}, \hat{\zeta})^2 - 1)} } = \\
&=\frac{\frac{2r}{\kappa((\theta_{||}, \hat{\zeta})^2 - 1)} \alpha}{ \frac{r|k_{||}-\kappa \theta_{||}|}{\kappa((\theta_{||}, \hat{\zeta})^2 - 1)} +\sqrt{\frac{r^2|k_{||}-\kappa \theta_{||}|^2}{\kappa^2((\theta_{||}, \hat{\zeta})^2 - 1)^2}+\frac{2r}{\kappa((\theta_{||}, \hat{\zeta})^2 - 1)} \alpha}} = \frac{2 \alpha}{ |k_{||}-\kappa \theta_{||}| +\sqrt{|k_{||}-\kappa \theta_{||}|^2+\frac{2 \kappa}{r} ((\theta_{||}, \hat{\zeta})^2 - 1) \alpha}}.
\end{aligned}
\end{align}
For simplicity, we consider \eqref{eq:4.14} assuming that $\alpha<0,$ and the square root is positive. 
Note that the argument of square root is positive for $\alpha<0$ in view of \eqref{eq:4.12}. 
Then, from \eqref{eq:4.14} it follows that
\begin{align}\label{eq:4.15}
    |\beta| \leq  \sqrt{ \frac{2\alpha r}{\kappa((\theta_{||}, \hat{\zeta})^2 - 1) }}.
\end{align}

Formula \eqref{eq:2.19} follows from \eqref{eq:2.14} and \eqref{eq:4.15}.

Note also that 
\begin{align}\label{eq:4.16}
   (\theta_{||}, \hat{\zeta})^2 - 1<  \delta^2 - 1, \text{ for } \theta \in \mathbb{S}^+_{\omega, \delta}.
\end{align}

Formula \eqref{eq:2.20} follows from  formula for $D$ in \eqref{eq:2.4}, formulas \eqref{eq:4.10.0}, \eqref{eq:4.11},  \eqref{eq:4.15}, and \eqref{eq:4.16}.

\subsection{Proof of Proposition \ref{prop:2.3}}

Proposition \ref{prop:2.3} follows from Theorem \ref{thm:1} and formulas \eqref{eq:2.9.0} and  \eqref{eq:2.9}. In details, for fixed $\kappa,$ this proof is as follows.

Due to \eqref{eq:2.9.0} and  \eqref{eq:2.9}, we have that:
\begin{align}
    &D = D(x, \zeta, k) = 2i \sin\left( \alpha + \frac{\alpha^2}{\varepsilon^2|x| } {\cal O}\left(1 \right)\right), \text{ as } \varepsilon |x|/\alpha \to +\infty, \text{ if } \hat{x} = \theta \in \mathbb{S}^+_{\omega} \setminus {\cal N}_{\omega, k, \varepsilon}, \\
    &D = 2i  \left(1 + {\cal O}\left( \frac{\alpha}{\varepsilon^2 |x|}\right)\right)\sin \alpha, \text{ as } \varepsilon |x|/ \alpha \to + \infty \, \text{ and } \varepsilon^2 |x|/\alpha \to +\infty, \quad  \text{ if }  \hat{x} = \theta \in \mathbb{S}^+_{\omega} \setminus {\cal N}_{\omega, k, \varepsilon}. \label{eq:4.9}
\end{align}
Theorem \ref{thm:1} and formula \eqref{eq:4.9} imply that 
\begin{align}
\begin{aligned}
 &f_1(\theta) = f_{1, 1}(\theta) + \frac{1}{2i \sin \alpha} \left(1 + {\cal O}\left( \frac{\alpha}{\varepsilon^2|x|}\right)\right)\left({\cal O}(\frac{1}{|x|^{\sigma}})+{\cal O}(\frac{\alpha}{\varepsilon|x|})\right) = \\
 &=f_{1, 1}(\theta) +   \frac{1}{2i \sin \alpha}{\cal O}\left( \frac{1}{|x|^{\sigma}}\right) + \frac{1}{2i \sin \alpha} \left({\cal O}(\frac{1}{\varepsilon|x|})+{\cal O} \left( \frac{\alpha}{\varepsilon^2 |x|^{1+\sigma}}\right)+{\cal O} \left( \frac{\alpha^2}{\varepsilon^3 |x|^2}\right)\right) , \\
 & \text{ as } \varepsilon |x|/ \alpha \to + \infty \, \text{ and } \varepsilon^2 |x|/\alpha \to +\infty, \quad  \text{ if }  \hat{x} = \theta \in \mathbb{S}^+_{\omega} \setminus {\cal N}_{\omega, k, \varepsilon}.
\end{aligned} \label{eq:4.10}
\end{align}

Note that, for fixed $\kappa,$ the property that $\varepsilon|x|/\alpha \to +\infty$ follows from the property that $\varepsilon^2 |x|/\alpha \to +\infty,$ since $\varepsilon \in (0, 2\kappa).$

Formula \eqref{eq:2.13} follows from \eqref{eq:4.10} at fixed $\alpha$ and $\varepsilon.$

Proposition \ref{prop:2.3} is proved.

\subsection{Proof of Proposition \ref{prop:2.5}}

Due to \eqref{eq:S+omegadelta}, \eqref{eq:2.19}, we have that 
\begin{align}\label{eq:4.20}
    |\zeta| \leq\sqrt{\frac{2|\alpha| r}{\kappa (1-\delta^2)}}, \quad \text{for } \theta \in \mathbb{S}^+_{\omega, \delta}.
\end{align}

Due to \eqref{eq:2.20} and the formula 
$\sin(\alpha+\beta) = \sin \alpha \cos \beta + \cos \alpha \sin \beta,$ we have that 
\begin{align}\label{eq:4.21}
    D =  2i  \left(1+ {\cal O}\left(\frac{1}{r^{1/2}}\right)\right)\sin \alpha, \quad \text{ as } r\to +\infty,
\end{align}
uniformly in $\theta \in \mathbb{S}^+_{\omega, \delta}.$

Theorem \ref{thm:1} and formulas \eqref{eq:4.20} and  \eqref{eq:4.21} imply that
\begin{align}\label{eq:4.22}
&f_1(\theta) = f_{1, 1}(\theta) + \frac{1}{2i \sin \alpha} \left[ 1 + {\cal O}\left( \frac{ 1}{r^{1/2}}\right) \right]\left[{\cal O}\left(\frac{1}{r^{\sigma}}\right)+{\cal O}\left(\frac{|\zeta|}{r}\right)\right] = f_{1, 1}(\theta) + \frac{1}{2i \sin \alpha} {\cal O}\left( \frac{1}{r^{1/2}} \right),    
\end{align}
as $r \to +\infty,$ uniformly in $\theta \in \mathbb{S}^+_{\omega, \delta}, $ where $x = r \theta.$

Formula \eqref{eq:2.14.1} follows from \eqref{eq:4.22}.

Proposition \ref{prop:2.5} is proved.

\section{Proof of Proposition \ref{prop:2.2} }\label{sec:5}

Due to \eqref{eq:3.2.0} for $d =2,$ we have that
\begin{align}\label{eq:5.1}
    a(x, k) = e^{-ikx+i\kappa |x|} f_1 + e^{ikx-i\kappa |x|} \overline{f_1} + \frac{1}{|x|^{1/2}} f_1 \overline{f_1} + {\cal O}\left(\frac{1}{|x|}\right), \quad \text{as } |x|\to +\infty,
\end{align}
where $f_1 = f_1(\theta),$ $\theta = \hat{x} = x/|x|, $ $x\in \mathbb{R}^2,$ uniformly in $\theta \in {\mathbb S}^{1}.$

In view of \eqref{eq:5.1}, formulas \eqref{eq:3.5}, \eqref{eq:3.6.0} admit the following more precise version, for $d=2$:
\begin{align}
&\begin{aligned}\label{eq:4.1}
&e^{-ikx}e^{i\kappa |x|} f_1(\theta) + e^{ikx}e^{-i\kappa |x|} \overline{f_1}(\theta) = a(x, k) - \frac{1}{|x|^{1/2}}f_1(\theta) \overline{f_1}(\theta) + {\cal O}\left(\frac{1}{|x|}\right), 
\end{aligned}\\
&\begin{aligned}\label{eq:4.2}
&e^{-iky}e^{i\kappa |y|} f_1(\theta) + e^{iky}e^{-i\kappa |y|}\overline{f_1}(\theta) = a(y, k) - \frac{1}{|y|^{1/2}}f_1(\hat{y}) \overline{f_1}(\hat{y})+ {\cal O}\left(\frac{1}{|x|}\right)+ {\cal O}\left(\frac{|\zeta|}{|x|}\right), 
\end{aligned}
\end{align}
 where $|x| \to +\infty$ and  $|x|/|\zeta|\to +\infty,$ uniformly in $\theta\in {\mathbb S}^+_{\omega};$  here, $x,$ $y$ are defined as in \eqref{eq:2.2}. 
In addition:
\begin{align}
&\frac{1}{|y|^{1/2}}f_1(\hat{y}) \overline{f_1}(\hat{y}) = \frac{1}{|x|^{1/2}}\left(1+{\cal O}\left(\frac{|\zeta|}{|x|}\right)\right)f_1(\hat{y}) \overline{f_1}(\hat{y}) =  \frac{1}{|x|^{1/2}}f_1(\hat{y}) \overline{f_1}(\hat{y}) + {\cal O}\left(\frac{|\zeta|}{|x|^{3/2}}\right),   
\end{align}
 in view of \eqref{eq:3.9.1};
\begin{align}\label{eq:4.3}
    f_1(\hat{y}) \overline{f_1}(\hat{y}) = \left(f_{1}(\theta) +{\cal O}\left(\frac{|\zeta|}{ |x|} \right) \right)\overline{\left(f_{1}(\theta) +{\cal O}\left(\frac{|\zeta|}{ |x|} \right) \right)} =   f_{1}(\theta) \overline{f_{1}}(\theta) + {\cal O}\left(\frac{|\zeta|}{|x|}\right),
\end{align}
 in view of \eqref{eq:3.6}. Thus, \eqref{eq:4.2}-\eqref{eq:4.3} imply that
 \begin{align}\label{eq:5.6}
     &e^{-iky}e^{i\kappa |y|} f_1(\theta) + e^{iky}e^{-i\kappa |y|}\overline{f_1}(\theta) = a(y, k) - \frac{1}{|x|^{1/2}}f_1(\theta) \overline{f_1}(\theta)+ {\cal O}\left(\frac{1}{|x|}\right)+ {\cal O}\left(\frac{|\zeta|}{|x|}\right),
 \end{align}
 as $|x|\to +\infty,$ $|x|/|\zeta|\to +\infty,$ uniformly in $\theta \in \mathbb{S}^+_{\omega}.$

In addition, in view of \eqref{eq:2.4} for $d=2$, we have that
\begin{align}\label{eq:5.7}
&f_{1, 1}(\theta) = \frac{1}{D}{\cal O}\left( 1\right),\\
&\begin{aligned}
&    f_1(\theta) \overline{f_1}(\theta) = \left(f_{1, 1}(\theta) +\frac{1}{D}\left({\cal O}\left(\frac{1}{ |x|^{1/2}}\right)+{\cal O}\left(\frac{|\zeta|}{ |x|}\right)\right)\overline{\left(f_{1, 1}(\theta) +\frac{1}{D}\left({\cal O}\left(\frac{1}{ |x|^{1/2}}\right)+{\cal O}\left(\frac{|\zeta|}{ |x|}\right)\right)\right)} \right) = \\
&=f_{1, 1}(\theta) \overline{f_{1, 1}}(\theta) + \frac{1}{D^2} \left({\cal O}\left(\frac{1}{ |x|^{1/2}}\right)+{\cal O}\left(\frac{|\zeta|}{ |x|}\right) \right)+ \frac{1}{D^2}\left({\cal O}\left(\frac{1}{ |x|^{1/2}}\right)+{\cal O}\left(\frac{|\zeta|}{ |x|}\right) \right)^2 = \\
&= f_{1, 1}(\theta) \overline{f_{1, 1}}(\theta) + \frac{1}{D^2} \left({\cal O}\left(\frac{1}{ |x|^{1/2}}\right)+{\cal O}\left(\frac{|\zeta|}{ |x|}\right) \right),
\end{aligned}
\end{align}
as $|x|\to+\infty,$ $|x|/|\zeta|\to +\infty,$ uniformly in $\theta \in \mathbb{S}^+_{\omega}.$  
 
Substituting \eqref{eq:5.7} into \eqref{eq:4.1}, \eqref{eq:5.6}, we obtain the following system for $f_1$ and $\overline{f_1}$:
\begin{align}
&\begin{aligned}\label{eq:4.4}
&e^{-ikx}e^{i\kappa |x|} f_1(\theta) + e^{ikx}e^{-i\kappa |x|} \overline{f_1}(\theta) = a(x, k) - \frac{1}{|x|^{1/2}}f_{1, 1}(\theta) \overline{f_{1, 1}}(\theta) + \\
&+{\cal O}\left(\frac{1}{|x|}\right) +  \frac{1}{D^2} \left({\cal O}\left(\frac{1}{ |x|}\right)+{\cal O}\left(\frac{|\zeta|}{ |x|^{3/2}}\right) \right), 
\end{aligned}\\
&\begin{aligned}\label{eq:4.5}
&e^{-iky}e^{i\kappa |y|} f_1(\theta) + e^{iky}e^{-i\kappa |y|}\overline{f_1}(\theta) = a(y, k) - \frac{1}{|x|^{1/2}}f_{1, 1}(\theta) \overline{f_{1, 1}}(\theta)+ \\
&+{\cal O}\left(\frac{1}{|x|}\right)+{\cal O}\left(\frac{|\zeta|}{|x|} \right) 
+\frac{1}{D^2} \left({\cal O}\left(\frac{1}{ |x|}\right)+{\cal O}\left(\frac{|\zeta|}{ |x|^{3/2}}\right) \right), 
\end{aligned}
\end{align}
as $|x|\to+\infty,$ $|x|/|\zeta|\to +\infty,$ uniformly in $\theta \in \mathbb{S}^+_{\omega}.$

Solving the system \eqref{eq:4.4}, \eqref{eq:4.5}, we obtain that
\begin{align}
f_1(\theta) = f_{1, 1}(\theta) - \frac{1}{D}\left(e^{iky-i\kappa|y|} - e^{ikx - i\kappa|x|} \right) \frac{1}{|x|^{1/2}}f_{1, 1}(\theta) \overline{f_{1, 1}}(\theta) + \Sigma,    
\end{align}
where 
\begin{align}
\begin{aligned}
&\Sigma=  \frac{1}{D}\left({\cal O}\left(\frac{1}{|x|}\right)+{\cal O}\left(\frac{|\zeta|}{|x|}\right)\right) + \frac{1}{D^3} \left({\cal O}\left(\frac{1}{|x|}\right)+{\cal O}\left(\frac{|\zeta|}{|x|^{3/2}}\right) \right) = \\
&=\frac{1}{D}{\cal O}\left(\frac{|\zeta|}{|x|}\right)+ \frac{1}{D^3}\left({\cal O}\left(\frac{1}{|x|}\right)+{\cal O}\left(\frac{|\zeta|}{|x|^{3/2}}\right)\right),
\end{aligned}
\end{align}
as $|x|\to+\infty,$ $|x|/|\zeta|\to +\infty,$ uniformly in $\theta \in \mathbb{S}^+_{\omega}.$

 Proposition \ref{prop:2.2}  is proved.


\section{Numerical implementation}\label{sec:6}

In this section we illustrate numerically our two-point  phase recovering from holographic data on a single plane proceeding from Proposition \ref{prop:2.5} and formula \eqref{eq:2.25} for the case when $\hat{k} = k/|k| = \omega.$  

In our numerical example, we assume that:
\begin{align}\label{eq:6.1}
\kappa = 4, \quad k = (4, 0, 0), \quad X = X_{s, \omega}, \text{ where } s=100, \quad \omega = e_1 = (1, 0, 0),  
\end{align}

\begin{align}\label{eq:6.2}
    \psi_1(x) = c\frac{e^{i\kappa|x-x^0|}}{|x-x^0|}, \quad c=1, \quad x^0 = (0, x^0_2, 0), \quad x^0_2 = 2.5. 
\end{align} 

\begin{align}
    {\cal U} = \{x \in \mathbb{R}^3: |x|>r \}, \text{ where } 2.5<r<100.
\end{align}

The reconstruction of $\psi_1$ from holographic data is considered on the square $ G  \subset X,$ where
\begin{align}\label{eq:G_def}
    G = \{x\in X: |x_2|< h, \, |x_3|<h \}, \quad h = 20.
\end{align}


In our numerical example, $G$ is considered as an $N\times N$ square grid, where $N=100.$ 

Figure 1 shows the intensity (hologram) $I = |\psi_0+\psi_1|^2$ restricted to $G.$ 
 \begin{figure}[H]
    \centering
    \includegraphics[width=0.44\linewidth]{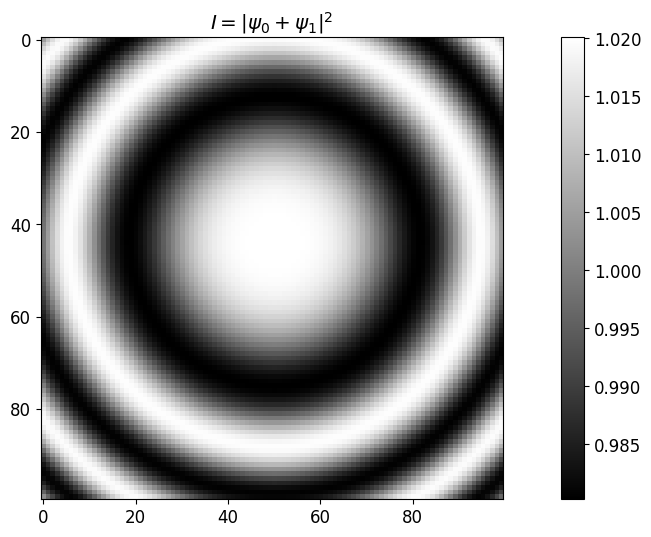}
    \caption{Intensity (hologram) $I = |\psi_0+\psi_1|^2.$}
\end{figure}

We reconstruct $\psi_1$ on $G$ as $\psi_{1, rec},$ where 
\begin{align}\label{eq:6.6}
&\psi_{1, rec}(x) = \frac{e^{i\kappa|x|}}{|x|^{\frac{d-1}{2}}}f_{1, 1}\left(\hat{x}\right), \quad x \in G,
\end{align}
 where $f_{1, 1}(\hat{x})$ is defined according to \eqref{eq:2.4} and \eqref{eq:2.2} with $\zeta = y -x$ defined proceeding from   \eqref{eq:2.14}, \eqref{eq:2.15}, \eqref{eq:2.18}. 
  In our numerical example, we use formulas \eqref{eq:2.14}, \eqref{eq:2.15}, \eqref{eq:2.18},  for $\alpha= -1/2.$ 
  In addition, in our numerical example, we have that $|\zeta|<\varepsilon = 15$ on $G.$ 

Note that in our numerical example, in view of $k$ and $\omega$ indicated in \eqref{eq:6.1}, we have that $k_{||} = 0$ in formulas \eqref{eq:2.14}, \eqref{eq:2.15}, \eqref{eq:2.18}, and the central point ${\mathfrak O}$ of $G,$ where $x_2 = x_3 =0,$ corresponds to $\theta_{||} = 0.$

Therefore, in our example, the central point ${\mathfrak O}$ corresponds to the singular  $\theta = k/\kappa$ arising in item (i) of Lemma \ref{lem:2.2}, and we  consider ${\mathfrak O}$  as a singular point on $G.$

\begin{figure}[H]
    \centering
    \includegraphics[width=0.9\linewidth]{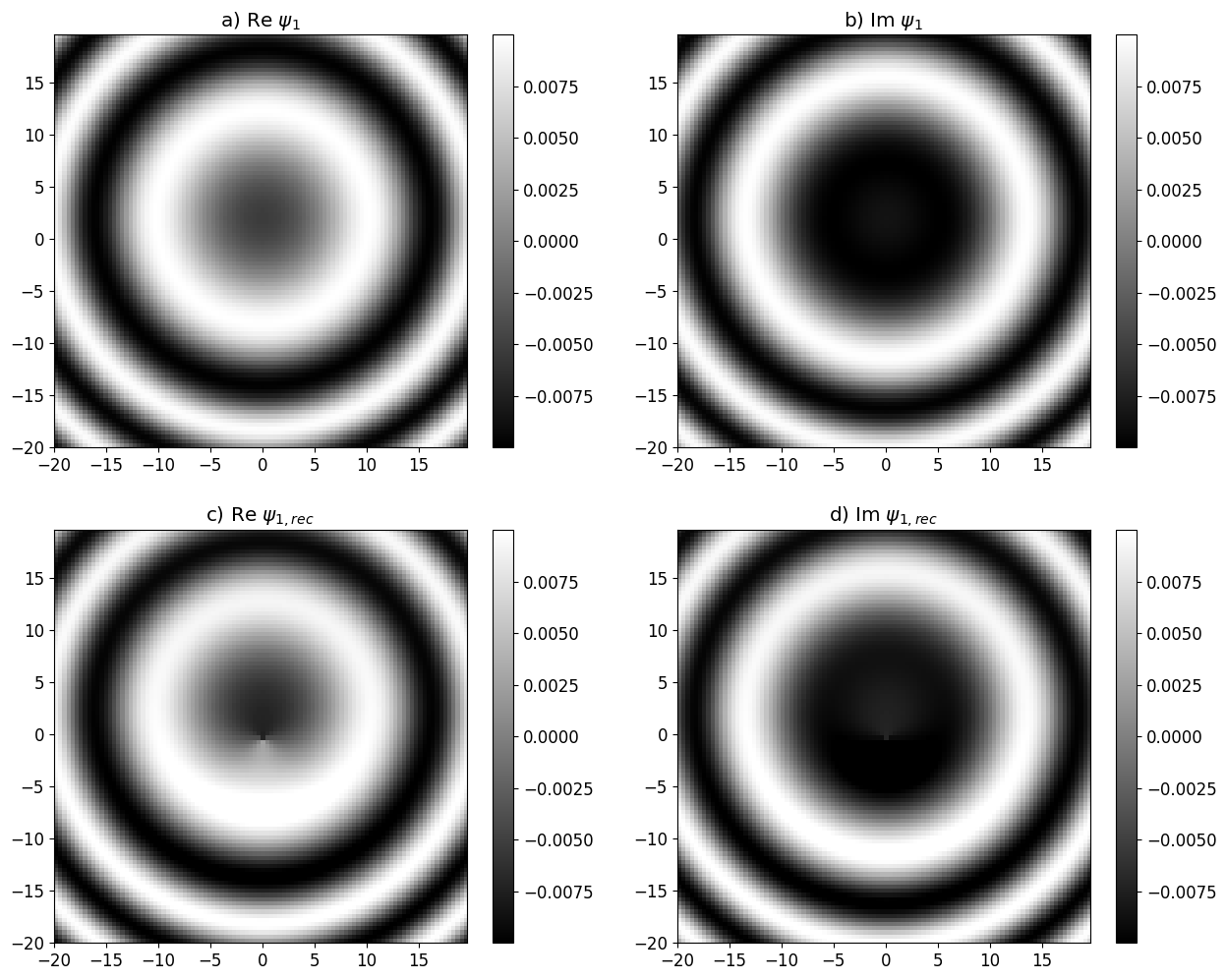}
    \caption{Exact $\psi_1$ and reconstructed  $\psi_{1, rec}$ on $G\subset X.$}
    \label{fig:enter-label-0}
\end{figure}

Figure \ref{fig:enter-label-0} shows $\psi_1$  and $\psi_{1, rec}$ on $G,$ where $\psi_{1, rec}$ is reconstructed from the intensity $|e^{ikx}+ \psi_1|^2$ on $G$ as described for formula \eqref{eq:6.6}. 

\begin{figure}[H]
    \centering
    \includegraphics[width=0.9\linewidth]{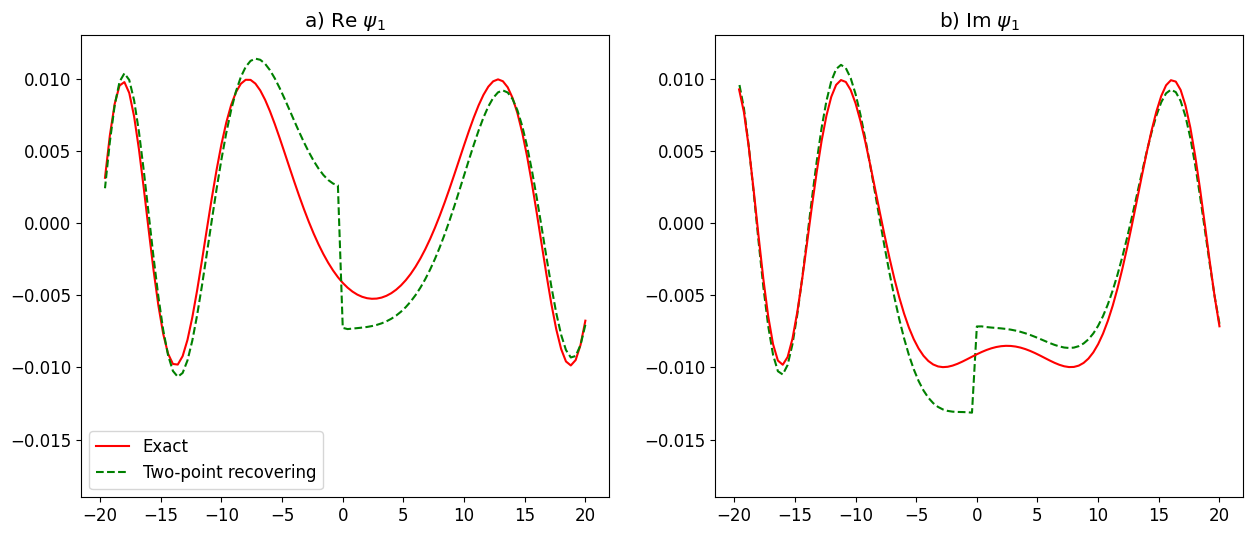}
    \caption{Exact (solid red) and recovered (dashed green) values of $Re \,\psi_1$ and $Im\,\psi_1$ along the central vertical line. }
    \label{fig:enter-label}
\end{figure}

 Figure \ref{fig:enter-label} shows  central vertical  profiles of $\psi_1$ and $\psi_{1, rec}$ on $G.$

Let 
\begin{align}
{\mathfrak E}(u_2, u_1, D) = \frac{\|u_2-u_1\|_{L_2(D)}}{\|u_1\|_{L_2(D)}},    
\end{align}
where $u_1,$ $u_2$ are some functions on $D \subseteq G,$ $u_1 \not\equiv 0$ on $D.$ 

In our example, we have that 
\begin{align}\label{eq:errors}
 {\mathfrak E}(\psi_{1, rec}, \psi_1,  G)   = 11.7\%, \quad {\mathfrak E}(\psi_{1, rec}, \psi_1, D)   = 29.7\%, \quad {\mathfrak E}(\psi_{1, rec}, \psi_1,  G\setminus D) = 10.2\%,    \quad D =  \{x\in G: |x_2|,|x_3|<2\}.   
\end{align}

Figures \ref{fig:enter-label-0}, \ref{fig:enter-label} and formula \eqref{eq:errors} show that our reconstruction $\psi_{1, rec}$  approximates $\psi_1$ with reasonable precision outside a neighbourhood 
of the singular point  ${\mathfrak O}.$ 






In addition, the quality of reconstruction $\psi_{1, rec}$ depends as follows on the  parameters $s, c, \kappa, x_2^0 $ in \eqref{eq:6.1}, \eqref{eq:6.2}:
\begin{itemize}\label{eq:list}
   \item  $s=5, $ ${\mathfrak E}= 25\%$; $s=10,$ ${\mathfrak E}= 16\%$; $s=100,$ ${\mathfrak E}= 11.7\%$; $s= 200, $ ${\mathfrak E}= 10.8\%;$
   \item   $c=0.1, {\mathfrak E} = 11.7\%;$ $c=1, {\mathfrak E} = 11.7\%;$   $c=10, {\mathfrak E} = 11.7\%;$ $c=20, {\mathfrak E} = 11.8\%;$
   \item  $\kappa = 1,{\mathfrak E} = 9.8\%;$ $   \kappa = 4, {\mathfrak E} = 11.7\%;$ $ \kappa =  16, {\mathfrak E} = 13.0\%;$ 
    \item  $x^0_2 = 0,  {\mathfrak E} = 0.13\%; $ $x^0_2 = 2.5,  {\mathfrak E} = 11.7\%; $ 
 $ x^0_2 = 5, {\mathfrak E} = 22.2\%.$ 
\end{itemize}
Here, ${\mathfrak E} = {\mathfrak E}(\psi_{1, rec}, \psi_1, G),$ where $\psi_1,$ $\psi_{1, rec}$ depend on $s, c, \kappa, x_2^0,$ and in each line above  we vary only one of these  parameters keeping other parameters as in \eqref{eq:6.1}, \eqref{eq:6.2}.

In addition to the reconstruction error ${\mathfrak E} = {\mathfrak E}(\psi_{1, rec}, \psi_1, D) $ we also consider
the discrepancy 
\begin{align}
    &\mathfrak{E}_{dis, D}  = \mathfrak{E}(I_{rec}-1, I-1, D), \quad I_{rec} = |\psi_0+\psi_{1, rec}|^2, \quad I = |\psi_0+\psi_1|^2,
\end{align}
where $I$ is our hologram restricted to $D.$ Here, we assume that $I \not\equiv  1$ on $D$ (since $\psi_1\equiv 0$ on $X$ if $I\equiv 1$ on $D,$ where $D$ is a non-empty open domain of $X,$ in view of \cite{N2023}). In our example we have that

\begin{align}\label{eq:discrepances}
 {\mathfrak E}_{dis, G}   = 7.2\cdot 10^{-3}, \quad {\mathfrak E}_{dis, D}   = 6.7\cdot 10^{-3}, \quad {\mathfrak E}_{dis, G\setminus D} = 7.2\cdot 10^{-3}, 
\end{align}
where $G,$ $D$ are as in \eqref{eq:G_def}, \eqref{eq:errors}.

Therefore,  the  discrepancy smallness  is not a proper criterion for good holographic reconstruction in the framework of our algorithm.

Our numerical examples also show that the reconstruction according to \eqref{eq:6.6} is not invariant with respect to parallel translations of the hologram in the plane $X.$  One can see that the singular point $\mathfrak{O}$ is related to the fixed coordinate system.
Therefore, probably, our formulas for finding $\psi_1$ on $X$ from $I$ on $X$ can be improved.

\section{Conclusions}\label{sec:conc}


In this work, we  contribute to holographic studies by presenting  a two-point approximate reconstructions of  an object beam $\psi_1$ on $X$ from the intensity (hologram) $I = |\psi_0+\psi_1|^2$ on $X,$ where $\psi_0$ is a reference beam and $X$ is a hyperplane in $\mathbb{R}^d,$ $d\geq 2$. We assume that $\psi_0$ and $\psi_1$ are modeled as  plane wave solution and radiation solution, respectively, for the Helmholtz equation \eqref{eq:schrod} in exterior region ${\cal U},$ and that $X$ is sufficiently far from the origin in $\mathbb{R}^d$. Our reconstructions are given in terms of the far-field pattern of $\psi_1$ (i.e., $f_1$ arising in formula \eqref{eq:1.5}) and are based on Theorem 2.1 and Corollary 2.1, and subsequent  results given in Section \ref{sec:main}. In these studies we proceed from the recent works \cite{NSh}, \cite{N2023}, \cite{NN}.

Our holographic reconstructions are given by explicit (elementary) formulas for approximate recovering  complex number  $\psi_1(x)$ for each  $x\in X$ from two real numbers $I(x)$ and $I(y),$ where $y \in X$ satisfies  some natural restrictions, implying that the determinant $D \neq 0$ and that the error term is relatively small in Theorem 2.1. These restrictions remain a lot of freedom for choosing $y.$ To our knowledge, the aforementioned  elementary formulas are much simpler than 
related formulas given in the literature for approximate finding $\psi_1$ on $X$ from the intensity  $I$ on $X.$ Note that in the present work we do not use any approximation for the Helmholtz equation, whereas approximations for the Helmholtz equation (e.g., paraxial approximation) are often used in holographic reconstructions. 

 We successfully implemented and illustrated numerically our two-point phase recovering  from holographic data on a single plane in dimension $d=3;$ see Section \ref{sec:6}. In addition, this implementation was realized for the case when the second point $y = x +\zeta$ was taken according to formulas \eqref{eq:2.14}, \eqref{eq:2.15}, \eqref{eq:2.18} with $\alpha = -1/2$. Although in our numerical example the reconstruction  looks good in general, see Figures 2, 3, some errors are still considerable,   at least, at some regions of the plane $X,$ see formula \eqref{eq:errors}. This seems to be related to the observation  (mentioned at the end of Section \ref{sec:6}) that the reconstruction according to \eqref{eq:6.6} is not invariant with respect to parallel translations of the hologram in the plane $X.$  Therefore, theoretical and numerical  studies of the present work should be continued.

\section{Acknowledgements}

V. N. Sivkin is supported by Russian Science Foundation, grant № 25-71-00116, https://rscf.ru/project/25-71-00116/

\bibliographystyle{alpha}

\vskip 3mm

Roman G. Novikov

CMAP, CNRS, \'Ecole polytechnique, Institut Polytechnique de Paris, 91128 Palaiseau, France 

\& IEPT RAS, 117997 Moscow, Russia

E-mail: roman.novikov@polytechnique.edu
\vskip 3mm

Vladimir N. Sivkin 

HSE University, Moscow, Russia

E-mail: sivkin96@yandex.ru

\end{document}